# SYMMETRY ANALYSIS OF THE CHARGED SQUASHED KALUZA-KLEIN BLACK HOLE METRIC

ROHOLLAH BAKHSHANDEH–CHAMAZKOTI

ABSTRACT. In this paper, a complete analysis of symmetries and conservation laws for the charged squashed Kaluza–Klein black hole spacetime in a Riemannian space is discussed. First, a comprehensive group analysis of the underlying space-time metric using Lie point symmetries are presented and then it the $n$-dimensional optimal system of this space-time metric, for $n = 1, \ldots, 4$, are computed. It is shown that there is not any $n$-dimensional optimal system of Lie symmetry subalgebra associated to the system of geodesic for $n \geq 5$. Then the point symmetries of the one parameter Lie groups of transformations that leave invariant the action integral corresponding to the Lagrangian that means Noether symmetries are found and then the conservation laws associated to the system of geodesic equations are calculated via Noether's theorem.

## Contents

## 1. Introduction

Symmetries of an equation are closely related to conservation laws. Noether's theorem [1, 2] provides a method for finding conservation laws of differential equations arising from a known Lagrangian and having a known Lie symmetry. A systematic and well known way of determining conservation laws for systems of Euler-Lagrange equations once their Noether symmetries are known is via Noether theorem. This theorem relies on the availability of a Lagrangian and the corresponding Noether symmetries which leave invariant the action integral, [3, 4].

Since the geodesic equations follow from the variation of the geodesic Lagrangian defined by the metric and due to the fact that the Noether symmetries are a subgroup of the Lie group of Lie symmetries of these equations, one should expect a relation of the Noether symmetries of this Lagrangian with the projective collineations of the metric or with its degenerates, [5, 6]. Recent work in this direction has been by Bokhari et al. [7, 8, 9] in which the relation of the Noether symmetries with the Killing vectors of some special spacetimes it is discussed.

In [5], Tsamparlis and Paliathanasis have computed the Lie point symmetries and the Noether symmetries explicitly together with the corresponding linear and quadratic first integrals for the Schwarzschild spacetime and the Friedman Robertson Walker (FRW) spacetime. These authors, in another paper [10], have proved a theorem which relates the Lie symmetries of the geodesic equations in a Riemannian space with the collineations of the metric. They applied the results to Einstein spaces and spaces of constant curvature.

Now in present paper, I try to find Lie point symmetry group and Noether symmetries for the charged squashed Kaluza–Klein black hole space-time is considered as the metric [11, 12]

$$(1.1) \qquad ds^2 = -f(r)dt^2 + \frac{k^2(r)}{f(r)}dr^2 + \frac{r^2}{4}[k(r)(\sigma_1^2 + \sigma_2^2) + \sigma_3^2],$$

where

$$\sigma_1 = \cos\psi\, d\theta + \sin\psi\, \sin\theta\, d\phi,$$
$$\sigma_2 = -\sin\psi\, d\theta + \cos\psi\, \sin\theta\, d\phi,$$
$$\sigma_3 = d\psi + \cos\theta\, d\phi.$$

The coordinates run in range, $0 \leq \theta < \pi$, $0 \leq \phi < 2\pi$, $0 \leq \psi < 2\pi$ and $0 < r < r_\infty$. Moreover, the functions in the metric define as

$$(1.2) \qquad f(r) = 1 - \frac{2M}{r^2} + \frac{q^2}{r^4}, \quad k(r) = \frac{f(r_\infty) r_\infty^4}{(r^2 - r_\infty^2)^2}.$$

Here $M$ and $q$ are the mass and charge of the black hole respectively.

This paper is organized as follows. In Section 2, a brief introduction of Lie symmetry for a geodesic equations is given and then I prove theorem 1 that is showed the charged squashed Kaluza-Klein black hole metric (1.1) has seventh order dimensional Lie point symmetry algebra. In Section 3, $n$th order dimensional optimal system for $n = 1, 2, 3, 4$ is obtained while there is not any $n$-dimensional optimal system of Lie subalgebra of $\mathfrak{g}$ associated to the system of geodesic for $n \geq 5$. In Section 4, the point generators (Noether symmetries) of the one parameter Lie groups of transformations that leave invariant the action integral corresponding to the Lagrangian are found. In Section 5, the conservation laws associated to the system of geodesic equations (2.11) are computed via Noether's theorem.

## 2. Lie point symmetries

Suppose $(M, g)$ is a Riemannian manifold of dimension $n$. In a local space-time coordinate like that $\mathbf{x}(s) = (x^1(s), \ldots, x^n(s))$, the geodesic equations form a system of $n$ nonlinear, second order ordinary differential equations

$$(2.1) \qquad \ddot{x}^i(t) + \Gamma^i_{jk} \dot{x}_j(t) \dot{x}_k(t) = 0, \qquad i = 1, 2, \ldots, n.$$

where $\Gamma^i_{jk}$ are the Christoffel symbols and "$\cdot$" represents derivative with respect to arc length $s$. Consider an second order system of $n$ ordinary differential equations (2.1) with following form

$$(2.2) \qquad R_i(s, \mathbf{x}(s), \mathbf{x}^{(1)}(s), \mathbf{x}^{(2)}(s)) = 0, \quad i = 1, \ldots, n.$$

where $\mathbf{x}^{(j)}(s)$, $j = 1, 2$, is the $j$th order derivative with respect to $s$. We consider a one-parameter Lie group of transformations acting on $(s, \mathbf{x})$-space as follows

$$(2.3) \qquad \bar{s} = s + \epsilon \xi(s, \mathbf{x}(s)) + O(\epsilon^2), \qquad \bar{x}^\alpha(s) = x^\alpha(s) + \epsilon \phi^\alpha(s, \mathbf{x}(s)) + O(\epsilon^2),$$

where $\alpha = 1, \ldots, n$. The infinitesimal generator $\mathbf{X}$ associated with the group of transformations (2.3) is

$$(2.4) \qquad \mathbf{X} = \xi(s, \mathbf{x}) \frac{\partial}{\partial s} + \eta^\alpha(s, \mathbf{x}) \frac{\partial}{\partial x^\alpha},$$

the second order prolongation of $\mathbf{X}$ is given by

$$(2.5) \qquad \mathbf{X}^{[2]} = \mathbf{X} + \eta^\alpha_{,(1)}(s, \mathbf{x}, \mathbf{x}^{(1)}) \frac{\partial}{\partial x^\alpha_{,(1)}} + \eta^\alpha_{,(2)}(s, \mathbf{x}, \ldots, \mathbf{x}^{(2)}) \frac{\partial}{\partial x^\alpha_{,(2)}},$$

in which the prolongation coefficients are

$$(2.6) \qquad \eta^\alpha_{,(1)} = D\eta^\alpha_{,(0)} - x^\alpha_{,(1)} D\xi, \qquad \eta^\alpha_{,(2)} = D\eta^\alpha_{,(1)} - x^\alpha_{,(2)} D\xi$$

where $\eta^\alpha_{,(0)} = \eta^\alpha(s, \mathbf{x})$ and $D$ is the total derivative operator.

The invariance of the (2.2) system under the one-parameter Lie group of transformations (2.3) leads to the invariance criterions [2, 13, 14, 15]. So $\mathbf{X}$ is a point symmetry generator of (2.2) if and only if

$$(2.7) \qquad \mathbf{X}^{[2]} R_i \big|_{R_i=0} = 0.$$

Using (2.7) we find a system of partial differential equations that is called the *determining equations*. Solving these determining equations leads the Lie point symmetries of (2.2).

**Theorem 1.** *The Lie algebra $\mathfrak{g}$ of Lie point symmetries associated to the charged squashed Kaluza-Klein black hole metric (1.1) has the following basis*

$$(2.8) \qquad \begin{aligned} &\mathbf{X}_1 = \frac{\partial}{\partial s}, \qquad \mathbf{X}_2 = s\frac{\partial}{\partial s}, \qquad \mathbf{X}_3 = \frac{\partial}{\partial t}, \qquad \mathbf{X}_4 = \frac{\partial}{\partial \varphi}, \qquad \mathbf{X}_5 = \frac{\partial}{\partial \psi}, \\ &\mathbf{X}_6 = \sin\varphi \frac{\partial}{\partial \theta} + \cos\varphi \cot\theta \frac{\partial}{\partial \varphi} - \cos\varphi \csc\theta \frac{\partial}{\partial \psi}, \\ &\mathbf{X}_7 = \cos\varphi \frac{\partial}{\partial \theta} - \sin\varphi \cot\theta \frac{\partial}{\partial \varphi} + \sin\varphi \csc\theta \frac{\partial}{\partial \psi}, \end{aligned}$$

*whose nonzero commutators are*

$$(2.9) \qquad [\mathbf{X}_1, \mathbf{X}_2] = \mathbf{X}_1, \quad [\mathbf{X}_4, \mathbf{X}_6] = \mathbf{X}_7, \quad [\mathbf{X}_4, \mathbf{X}_7] = -\mathbf{X}_6, \quad [\mathbf{X}_6, \mathbf{X}_7] = \mathbf{X}_4.$$

**Proof:** Consider the charged squashed Kaluza–Klein black hole space-time metric (1.1). The nonzero components of the Christoffel symbols for this metric are

$$(2.10) \qquad \begin{aligned} &\Gamma^t_{tr} = -\Gamma^r_{tt} = \frac{f'(r)}{2}, \quad \Gamma^r_{rr} = \frac{k(r)}{2f^2(r)}\left[2f(r)k'(r) - k(r)f'(r)\right], \\ &\Gamma^r_{\theta\theta} = \Gamma^\theta_{r\theta} = \frac{r}{8}\left[2k(r) + rk'(r)\right], \\ &\Gamma^r_{\varphi\varphi} = \Gamma^\varphi_{r\varphi} = \frac{r}{8}\left[(2k(r) + rk'(r))\sin^2\theta + 2\cos^2\theta\right], \\ &\Gamma^r_{\varphi\psi} = \Gamma^\varphi_{r\psi} = \Gamma^\psi_{r\varphi} = \frac{r}{2}\cos\theta, \quad \Gamma^\theta_{\varphi\psi} = \Gamma^\varphi_{\theta\psi} = \Gamma^\psi_{\theta\psi} = -\frac{r^2}{4}\sin\theta, \quad \Gamma^r_{\psi\psi} = \Gamma^\psi_{r\psi} = \frac{r}{4}, \end{aligned}$$

where $f' = \dfrac{df}{dr}$, $k' = \dfrac{dk}{dr}$. Now substituting the (2.10) symbols in (2.1) we obtain the following system of geodesic equations for the metric (1.1):

$$(2.11) \begin{cases} R_1: \ \ddot{t} + f'(r)\dot{t}\dot{r} = 0, \\ R_2: \ \ddot{r} - \dfrac{f'(r)}{2}\dot{t}^2 + \dfrac{k(r)}{2f^2(r)}\left[2f(r)k'(r) - k(r)f'(r)\right]\dot{r}^2 + \dfrac{r}{8}\left[2k(r) + rk'(r)\right]\dot{\theta}^2 \\ \qquad + \dfrac{r}{8}\left[(2k(r) + rk'(r))\sin^2\theta + 2\cos^2\theta\right]\dot{\varphi} + r\cos\theta\dot{\varphi}\dot{\psi} + \dfrac{r}{4}\dot{\psi}^2 = 0, \\ R_3: \ \ddot{\theta} + \dfrac{r}{4}\left[2k(r) + rk'(r)\right]\dot{r}\dot{\theta} - \dfrac{r^2}{2}\sin\theta\dot{\varphi}\dot{\psi} = 0, \\ R_4: \ \ddot{\varphi} + \dfrac{r}{4}\left[(2k(r) + rk'(r))\sin^2\theta + 2\cos^2\theta\right]\dot{r}\dot{\varphi} + r\cos\theta\dot{r}\dot{\psi} - \dfrac{r^2}{2}\sin\theta\dot{\theta}\dot{\psi} = 0, \\ R_5: \ \ddot{\psi} + r\cos\theta\dot{r}\dot{\varphi} - \dfrac{r^2}{2}\sin\theta\dot{\theta}\dot{\psi} + \dfrac{r}{2}\dot{r}\dot{\psi} = 0. \end{cases}$$

Let the vector field $\mathbf{X}$ with form

$$(2.12) \qquad \mathbf{X} = \xi\dfrac{\partial}{\partial s} + \tau\dfrac{\partial}{\partial t} + \rho\dfrac{\partial}{\partial r} + \Theta\dfrac{\partial}{\partial \theta} + \Phi\dfrac{\partial}{\partial \varphi} + \Psi\dfrac{\partial}{\partial \psi}$$

is the infinitesimal generator $\mathbf{X}$ associated with the group of transformations

$$\bar{s} \mapsto s + \epsilon\xi(s, \mathbf{x}(s)),$$
$$\bar{t} \mapsto t + \epsilon\tau(s, \mathbf{x}(s)),$$
$$\bar{r} \mapsto r + \epsilon\rho(s, \mathbf{x}(s)),$$
$$\bar{\theta} \mapsto \theta + \epsilon\Theta(s, \mathbf{x}(s)),$$
$$\bar{\phi} \mapsto \phi + \epsilon\Phi(s, \mathbf{x}(s)),$$
$$\bar{\psi} \mapsto \psi + \epsilon\Psi(s, \mathbf{x}(s)).$$

The second order prolongation of $\mathbf{X}$ is given by the (2.5) and (2.6) formulas. With imposing invariance criterions (2.7) on the (2.11), we find infinitesimal generator (2.8) which are Lie point symmetries of the charged squashed Kaluza–Klein black hole space-time (1.1), with (2.9) commutators relations. $\square$

3. The optimal system

Let $G$ be a Lie group, with $\mathfrak{g}$ its Lie algebra. Each element $T \in G$ yields inner automorphism $T_a \longrightarrow TT_aT^{-1}$ of the group $G$. Every automorphism of the group $G$ induces an automorphism of $\mathfrak{g}$. The set of all these automorphism is a Lie group called *the adjoint group* $G^A$. The Lie algebra of $G^A$ is the adjoint algebra of $\mathfrak{g}$, defined as follows. Let us have two infinitesimal generators $X, Y \in L$. The linear mapping $\text{Ad}X(Y) : Y \longrightarrow [X, Y]$ is an automorphism of $\mathfrak{g}$, called *the inner derivation of* $\mathfrak{g}$. The set of all inner derivations $\text{ad}X(Y)(X, Y \in \mathfrak{g})$ together with the Lie bracket $[\text{Ad}X, \text{Ad}Y] = \text{Ad}[X, Y]$ is a Lie algebra $\mathfrak{g}^A$ called the *adjoint algebra of* $\mathfrak{g}$. Clearly $\mathfrak{g}^A$ is the Lie algebra of $G^A$. Two subalgebras in $\mathfrak{g}$ are *conjugate* if there is a transformation of $G^A$ which takes one subalgebra into the other. The collection of pairwise non-conjugate $s$-dimensional subalgebras is the optimal system of subalgebras of order $s$. The construction of the one-dimensional optimal system of subalgebras can be carried out by using a global matrix of the adjoint transformations as suggested by Ovsiannikov [15]. The latter problem, tends to determine a list (that is called an *optimal system*) of conjugacy inequivalent subalgebras with the property that any other subalgebra is equivalent to a unique member of the list under some element of the adjoint representation i.e. $\overline{\mathfrak{h}} \operatorname{Ad}(g) \mathfrak{h}$ for some g of a considered Lie group. Thus we will deal with the construction of the optimal system of subalgebras of $\mathfrak{g}$. The adjoint action is given by the Lie series

$$(3.1) \qquad \operatorname{Ad}(\exp(s\,\mathbf{X}_i))\mathbf{X}_j = \mathbf{X}_j - s\,[\mathbf{X}_i, \mathbf{X}_j] + \frac{s^2}{2}\,[\mathbf{X}_i, [\mathbf{X}_i, \mathbf{X}_j]] - \cdots,$$

where $s$ is a parameter and $i, j = 1, \cdots, n$. We can expect to simplify a given arbitrary element,

$$(3.2) \qquad \mathbf{X} = \sum_{i=1}^{7} a_i \mathbf{X}_i,$$

of the geodesics system Lie algebra $\mathfrak{g}$. Note that the elements of $\mathfrak{g}$ can be represented by vectors $(a_1, \ldots, a_7) \in \mathbb{R}^7$ since each of them can be written in the form (3.2) for some constants $a_1, \ldots, a_7$. Hence, the adjoint action can be regarded as (in fact is) a group of linear transformations of the vectors $(a_1, \ldots, a_7)$.

**Theorem 2.** *An optimal system of one-dimensional Lie subalgebras associated to the system of geodesic is generated by*

(1)  $A_1^1 = \langle \mathbf{X}_1 + a\mathbf{X}_3 + b\mathbf{X}_5 + c\mathbf{X}_7 \rangle$ \qquad (4)  $A_1^4 = \langle \mathbf{X}_4 + a\mathbf{X}_2 + b\mathbf{X}_3 + c\mathbf{X}_5 + d\mathbf{X}_6 \rangle$

(2)  $A_1^2 = \langle \mathbf{X}_1 + a\mathbf{X}_3 + b\mathbf{X}_5 + c\mathbf{X}_6 \rangle$ \qquad (5)  $A_1^5 = \langle \mathbf{X}_7 + a\mathbf{X}_2 + b\mathbf{X}_3 + c\mathbf{X}_5 \rangle$

(3)  $A_1^3 = \langle \mathbf{X}_1 + a\mathbf{X}_3 + b\mathbf{X}_5 + c\mathbf{X}_4 \rangle$ \qquad (6)  $A_1^6 = \langle a\mathbf{X}_2 + b\mathbf{X}_3 + c\mathbf{X}_5 \rangle$

*where $a, b, c, d \in \mathbb{R}$ are arbitrary constants.*

*Proof.* The function $F_i^s : \mathfrak{g} \to \mathfrak{g}$ defined by $\mathbf{X} \mapsto \operatorname{Ad}(\exp(s_i \mathbf{X}_i).\mathbf{X})$ is a linear map, for $i = 1, \cdots, 7$. The matrix $M_i^s$ of $F_i^s$ with respect to basis $\{\mathbf{X}_1, \cdots, \mathbf{X}_7\}$ is

$$M_1^{s_1} = \begin{bmatrix} 1 & 0 & 0_{15} \\ -s_1 & 1 & 0_{15} \\ 0_{51} & 0_{51} & \mathbf{I}_5 \end{bmatrix}, \qquad M_2^{s_2} = \begin{bmatrix} \exp(s_2) & 0_{16} \\ 0_{61} & \mathbf{I}_{55} \end{bmatrix}, \qquad M_3^{s_3} = \mathbf{I}_7,$$

$$M_4^{s_4} = \begin{bmatrix} \mathbf{I}_5 & 0_{51} & 0_{51} \\ 0_{15} & \cos(s_4) & -\sin(s_4) \\ 0_{15} & -\sin(s_4) & \cos(s_4) \end{bmatrix}, \qquad M_5^{s_5} = \mathbf{I}_7,$$

$$M_6^{s_6} = \begin{bmatrix} \mathbf{I}_3 & 0_{31} & 0_{32} & 0_{31} \\ 0_{13} & \cos(s_6) & 0_{12} & -\sin(s_6) \\ 0_{23} & 0_{21} & \mathbf{I}_2 & 0 \\ 0_{113} & -\sin(s_6) & 0_{12} & \cos(s_6) \end{bmatrix},$$

$$M_7^{s_7} = \begin{bmatrix} 1 & 0 & 0 & 0 & 0 & 0 & 0 \\ -s_1 & 1 & 0 & 0 & 0 & 0 & 0 \\ 0 & 0 & 1 & 0 & 0 & 0 & 0 \\ 0 & 0 & 0 & \cos(s_7) & 0 & \sin(s_7) & 0 \\ 0 & 0 & 0 & 0 & 1 & 0 & 0 \\ 0 & 0 & 0 & \sin(s_7) & 0 & \cos(s_7) & 0 \\ 0 & 0 & 0 & 0 & 0 & 0 & 1 \end{bmatrix},$$

where $\mathbf{I}_n$ is $n \times n$ identity matrix and $0_{ij}$ is $i \times j$ zero matrix. Let $\mathbf{X} = \sum_{i=1}^{7} a_i \mathbf{X}_i$ then it is seen that

$$F_7^{s_7} \circ F_6^{s_6} \circ \cdots \circ F_1^{s_1} : \mathbf{X} \mapsto [\exp(s_2)a_1]\mathbf{X}_1 + [-s_1 \exp(s_2)a_1 + a_2]\mathbf{X}_2 + a_3\mathbf{X}_3$$
$$+[\cos(s_6)\cos(s_7)a_4 - \cos(s_6)\sin(s_7)a_6 + \sin(s_6)a_7]\mathbf{X}_4 + a_5\mathbf{X}_5 +$$
$$[(\sin(s_4)\sin(s_6)\cos(s_7) + \cos(s_4)\sin(s_7))a_4 + (-\sin(s_4)\sin(s_6)\sin(s_7)$$
$$+ \cos(s_4)\cos(s_7))a_6 - \sin(s_4)\cos(s_6)a_7]\mathbf{X}_6 + [(-\cos(s_4)\sin(s_6)\cos(s_7)$$
$$- \sin(s_4)\sin(s_7))a_4 + (\cos(s_4)\sin(s_6)\sin(s_7) - \sin(s_4)\cos(s_7))a_6 + \cos(s_4)\cos(s_6)a_7]\mathbf{X}_7.$$

If $a_1, a_7 \neq 0$ then we can omit the coefficients of $\mathbf{X}_2, \mathbf{X}_4$ and $\mathbf{X}_6$ by setting $s_1 = \frac{a_2}{a_1}$, $s_2 = s_7 = 0$, $s_4 = \arctan(\frac{a_6}{a_7})$ and $s_6 = -\arctan(\frac{a_4}{a_7})$. Scaling $\mathbf{X}$, we can assume that $a_1 = 1$. So, $\mathbf{X}$ is reduced to the case (1).

If $a_1, a_6 \neq 0$ and $a_7 = 0$ then the coefficients of $\mathbf{X}_2, \mathbf{X}_4$ and $\mathbf{X}_7$ are vanished by setting $s_1 = \frac{a_2}{a_1}$, $s_2 = s_4 = s_6 = 0$, and $s_7 = \arctan(\frac{a_4}{a_6})$. Scaling $\mathbf{X}$, we can assume that $a_1 = 1$. So, $\mathbf{X}$ is reduced to the case (2).

If $a_1 \neq 0$ and $a_6 = a_7 = 0$ then one can vanish the coefficients of $\mathbf{X}_2, \mathbf{X}_6$ and $\mathbf{X}_7$ by setting $s_1 = \frac{a_2}{a_1}$, $s_2 = s_6 = s_7 = 0$. Scaling $\mathbf{X}$, we can assume that $a_1 = 1$. So, $\mathbf{X}$ is reduced to the case (3).

If $a_1 = 0$ and $a_4, a_6 \neq 0$ then one can vanish the coefficient of $\mathbf{X}_7$ by setting $s_6 = s_7 = 0$ and $s_4 = \arctan(\frac{a_7}{a_6})$. Scaling $\mathbf{X}$, we can assume that $a_4 = 1$. So, $\mathbf{X}$ is reduced to the case (4).

If $a_1 = a_4 = a_6 = 0$ and $a_7 \neq 0$ then one can vanish the coefficients of $\mathbf{X}_4$ and $\mathbf{X}_6$ by setting $s_6 = s_7 = 0$. Scaling $\mathbf{X}$, we can assume that $a_7 = 1$. So, $\mathbf{X}$ is reduced to the case (5).

Finally, $a_1 = a_4 = a_6 = a_7 = 0$ leads the case (6). □

**Theorem 3.** *Two-dimensional optimal systems of Lie subalgebras associated to the system of geodesic have following forms:*

(1) $B_2^1 = \langle a\mathbf{X}_2 + b\mathbf{X}_3 + c\mathbf{X}_5, \mathbf{X}_1 \rangle$,  (2) $B_2^2 = \langle a\mathbf{X}_2 + b\mathbf{X}_3 + c\mathbf{X}_5, \mathbf{X}_7 \rangle$,

(3) $B_2^2 = \langle a\mathbf{X}_2 + b\mathbf{X}_3 + c\mathbf{X}_5, \mathbf{X}_6 + \alpha\mathbf{X}_7 \rangle$,

*where $a, b, c, \alpha$ are arbitrary constants.*

*Proof.* Suppose that $\mathfrak{h} = \mathrm{Span}_{\mathbb{R}}\{\mathbf{X}, \mathbf{Y}\}$ be a two-dimensional Lie subalgebra of $\mathfrak{g}$ with $\mathbf{X} = a\mathbf{X}_2 + b\mathbf{X}_3 + c\mathbf{X}_5$ and $\mathbf{Y} = \sum_{i=1}^{7} \gamma_i \mathbf{X}_i$, and then $[\mathbf{X}, \mathbf{Y}] = m\mathbf{X} + n\mathbf{Y}$, where $m, n \in \mathbb{R}$. If $n \neq 0$ then we have $\gamma_4 = \gamma_6 = \gamma_7 = 0$ and since $\mathfrak{h}$ has dimension two therefore we find case (1). If $n = m = 0$ then we have $\gamma_1 = 0$ and since $\dim \mathfrak{h} = 2$ thus $\mathbf{Y} = \gamma_4 \mathbf{X}_4 + \gamma_6 \mathbf{X}_6 + \gamma_7 \mathbf{X}_7$ where $\gamma_1^2 + \gamma_2^2 + \gamma_3^2 \neq 0$. Now if $\gamma_7 \neq 0$ then we can make the coefficients of $\mathbf{X}_4$ and $\mathbf{X}_6$ vanish by setting $s_7 = 0$, $s_4 = \arctan(\frac{a_6}{a_7})$ and $s_6 = -\arctan(\frac{a_4}{a_7})$. By Scaling $\mathbf{Y}$, we may assume $\gamma_7 = 1$ and then $\mathbf{Y} = X_7$. So, $\mathbf{Y}$ is reduced to the case (2). If $\gamma_6 \neq 0$ then we may remove the coefficient of $\mathbf{X}_4$ by putting $s_4 = s_6 = 0$, and $s_7 = \arctan(\frac{a_4}{a_6})$. Scaling $\mathbf{Y}$, we may choose $\gamma_4 = 1$ and then $\mathbf{Y} = \mathbf{X}_6 + \alpha \mathbf{X}_7$. So, $\mathbf{Y}$ is reduced to the case (3). There is not any new cases because choosing $\mathbf{X}$ equal to one of vector fields of cases (1)-(5) in Theorem 2 doesn't lead to $\dim \mathfrak{h} = 2$. □

**Theorem 4.** *Three-dimensional optimal systems of Lie subalgebras associated to the system of geodesic have following forms:*

(1) $C_3^1 = \langle a\mathbf{X}_2 + b\mathbf{X}_3 + c\mathbf{X}_5, \mathbf{X}_7, \mathbf{X}_4 - \mathbf{X}_6 \rangle,$   (2) $C_3^2 = \langle a\mathbf{X}_2 + b\mathbf{X}_3 + c\mathbf{X}_5, \mathbf{X}_7, X_1 \rangle,$

*where $a, b, c$ are arbitrary constants.*

*Proof.* Let $\mathfrak{l} = \text{Span}_{\mathbb{R}}\{\mathbf{X}, \mathbf{Y}, \mathbf{Z}\}$ be a three dimensional Lie subalgebra of $\mathfrak{g}$ with $\mathbf{X} = a\mathbf{X}_2 + b\mathbf{X}_3 + c\mathbf{X}_5$, $\mathbf{Y} = \mathbf{X}_7$ and $\mathbf{Z} = \sum_{i=1}^{7} \sigma_i \mathbf{X}_i$, and then $[\mathbf{X}, \mathbf{Z}] = m_1\mathbf{X} + m_2\mathbf{Y} + m_3\mathbf{Z}$ and $[\mathbf{Y}, \mathbf{Z}] = n_1\mathbf{X} + n_2\mathbf{Y} + n_3\mathbf{Z}$, where $m_i, n_i \in \mathbb{R}$. If $n_3 \neq 0$ then we have $\sigma_1 = 0$ and $\sigma_4 = \sigma_6$ and since $\dim\mathfrak{l} = 3$ so we obtain $Z = \mathbf{X}_4 - \mathbf{X}_6$ and we find case (1). If $n_3 = 0$, $m_3 \neq 0$ then we have $\sigma_4 = \sigma_6 = 0$ and thus $\mathbf{Z} = \mathbf{X}_1$ because $\dim\mathfrak{l} = 3$ and this is a new case (2).

Assume that $\mathfrak{h}$ is generated by $\mathbf{X} = a\mathbf{X}_2 + b\mathbf{X}_3 + c\mathbf{X}_5$, $\mathbf{Y} = \mathbf{X}_6 + \alpha\mathbf{X}_7$ and $\mathbf{Z} = \sum_{i=1}^{7} \delta_i \mathbf{X}_i$. Then $[\mathbf{X}, \mathbf{Z}] = r_1\mathbf{X} + r_2\mathbf{Y} + r_3\mathbf{Z}$ and $[\mathbf{Y}, \mathbf{Z}] = s_1\mathbf{X} + s_2\mathbf{Y} + s_3\mathbf{Z}$, where $r_i, s_i \in \mathbb{R}$. This case leads to $\delta_1 = \delta_4 = 0$ and then there is not an algebra with 3 dimension. There is not any new cases $\dim\mathfrak{h} = 3$. $\square$

**Theorem 5.** *Four-dimensional optimal system of Lie subalgebra associated to the system of geodesic have following form:*

$$D_4^1 = \langle a\mathbf{X}_2 + b\mathbf{X}_3 + c\mathbf{X}_5, \mathbf{X}_7, \mathbf{X}_4 - \mathbf{X}_6, \mathbf{X}_1 \rangle,$$

*where $a, b, c$ are arbitrary constants and $a \neq 0$.*

*Proof.* Assume that the $\mathfrak{p} = \text{Span}_{\mathbb{R}}\{\mathbf{X}, \mathbf{Y}, \mathbf{Z}, \mathbf{V}\}$ be a 4–dimensional Lie subalgebra of $\mathfrak{g}$ with $\mathbf{X} = a\mathbf{X}_2 + b\mathbf{X}_3 + c\mathbf{X}_5$, $\mathbf{Y} = \mathbf{X}_7$, $\mathbf{Z} = \mathbf{X}_4 - \mathbf{X}_6$ and $\mathbf{V} = \xi_1\mathbf{X}_4$, and then $[\mathbf{X}, \mathbf{V}] = a_1\mathbf{X} + a_2\mathbf{Y} + a_3\mathbf{Z} + a_4\mathbf{V}$, $[\mathbf{Y}, \mathbf{V}] = b_1\mathbf{X} + b_2\mathbf{Y} + b_3\mathbf{Z} + b_4\mathbf{V}$, and $[\mathbf{Z}, \mathbf{V}] = c_1\mathbf{X} + c_2\mathbf{Y} + c_3\mathbf{Z} + c_4\mathbf{V}$ where $a_i, b_i, c_i \in \mathbb{R}$. If $a_4 = 0$ then we have $\xi_1 = 0$ and this leads a contradiction with $\dim\mathfrak{l} = 4$ so we take $a_4 \neq 0$ and we have $a_4 = -a$ then $\mathbf{V} = \mathbf{X}_1$ in provided by $a \neq 0$. There is not any new cases $\dim\mathfrak{p} = 4$. $\square$

**Corollary 1.** *There is not any n-dimensional optimal system of Lie subalgebra of $\mathfrak{g}$ associated to the system of geodesic for $n \geq 5$.*

## 4. Noether symmetries

Consider an second order system of $m$ ordinary differential equations (2.1) with following form

(4.1) $$E_i(s, \mathbf{x}(s), \mathbf{x}^{(1)}(s), \mathbf{x}^{(2)}(s)) = 0, \quad i = 1, \ldots, n.$$

Let $\mathbf{Y}$ be a vector field defined on a real parameter fibre bundle over the manifold [14]

(4.2) $$\mathbf{Y} = \overline{\xi}(s, x^\mu)\frac{\partial}{\partial s} + \overline{\eta}^\nu(s, x^\mu)\frac{\partial}{\partial x^\nu},$$

where $\mu, \nu = 1, 2, 3, 4, 5$. The first prolongation of the above vector field is:

(4.3) $$\mathbf{Y}^{[1]} = \mathbf{Y} + \left(\overline{\eta}_{,s}^\nu + \overline{\eta}_{,\mu}^\nu \dot{x}^\mu - \overline{\xi}_{,s}\dot{x}^\nu - \overline{\xi}_{,s}\dot{x}^\mu\dot{x}^\nu\right)\frac{\partial}{\partial \dot{x}^\nu},$$

Then $\mathbf{X}$ is a *Noether point symmetry* of the Lagrangian

(4.4) $$\mathcal{L}(s, x^\mu, \dot{x}^\mu) = g_{\mu\nu}(x^\sigma)\dot{x}^\mu\dot{x}^\nu$$

if there exists a gauge function, $A(s, x^\mu)$, such that

(4.5) $$\mathbf{Y}^{[1]}\mathcal{L} + (D_s\xi)\mathcal{L} = D_sA,$$

where

$$(4.6) \quad D_s = \frac{\partial}{\partial s} + \dot{x}^\mu \frac{\partial}{\partial x^\mu}.$$

Since the geodesic equations are second order ordinary differential equations, one can take first order Lagrangian $\mathcal{L}(s, x_i, \dot{x}_i)$, where "dots" denotes differentiation with respect to the arc length parameter $s$, which leads a system of second ODEs

$$(4.7) \quad \ddot{x}^\mu = g(s, x^\mu, \dot{x}^\mu).$$

**Theorem 6.** *The Lie algebra of Noether symmetries associated to the charged squashed Kaluza-Klein black hole metric (1.1) has six dimension with following basis*

$$(4.8) \quad \begin{aligned} \mathbf{Y}_1 &= \frac{\partial}{\partial s}, \qquad \mathbf{Y}_2 = \frac{\partial}{\partial t}, \qquad \mathbf{Y}_3 = \frac{\partial}{\partial \varphi}, \qquad \mathbf{Y}_4 = \frac{\partial}{\partial \psi}, \\ \mathbf{Y}_5 &= \sin\varphi\, \frac{\partial}{\partial \theta} + \cos\varphi \cot\theta\, \frac{\partial}{\partial \varphi} - \cos\varphi \csc\theta\, \frac{\partial}{\partial \psi}, \\ \mathbf{Y}_6 &= \cos\varphi\, \frac{\partial}{\partial \theta} - \sin\varphi \cot\theta\, \frac{\partial}{\partial \varphi} + \sin\varphi \csc\theta\, \frac{\partial}{\partial \psi}. \end{aligned}$$

**Proof:** The Lagrangian for the (1.1) metric is

$$(4.9) \quad \mathcal{L} = -f(r)\dot{t}^2 + \frac{k^2(r)}{f(r)}\dot{r}^2 + \frac{r^2}{4}k(r)\dot{\theta}^2 + \frac{r^2}{4}\big[k(r)\sin^2\theta + \cos^2\theta\big]\dot{\varphi}^2 + \frac{r^2}{2}\cos\theta\,\dot{\varphi}\dot{\psi} + \frac{r^2}{4}\dot{\psi}^2,$$

where the "dot" represents derivative with respect to arc length $s$, and the functions $f(r)$ and $k(r)$ are introduced in (1.2). Then the Euler-Lagrange geodesic equations associated with the Lagrangian (4.9) are

$$(4.10) \quad \begin{cases} E_1: \ f(r)\ddot{t} + f'(r)\dot{t}\dot{r} = 0, \\ E_2: \ -\frac{2k^2(r)}{f(r)}\ddot{r} - f'(r)\dot{t}^2 + \frac{k(r)}{f^2(r)}[k(r)f'(r) - 2f(r)k'(r)]\dot{r}^2 + \frac{r}{4}[2k(r) + rk'(r)]\dot{\theta}^2 \\ \qquad + \frac{r}{4}[2k(r)\sin^2\theta + (2 + rk'(r))\sin^2\theta]\dot{\varphi}^2 + r\cos\theta\,\dot{\varphi}\dot{\psi} + \frac{r}{2}\dot{\psi}^2 = 0, \\ E_3: \ -\frac{r^2 k(r)}{2}\ddot{\theta} + \frac{r^2}{4}[(k(r)-1)\sin 2\theta]\dot{\varphi}^2 - \frac{r^2}{2}\sin\theta\,\dot{\varphi}\dot{\psi} - \frac{r}{2}[2k(r) + rk'(r)]\dot{r}\dot{\theta} = 0, \\ E_4: \ -\frac{r^2}{2}[k(r)\sin^2\theta + \cos^2\theta]\ddot{\varphi} - \frac{r^2}{2}\cos\theta\,\ddot{\psi} - \frac{r}{2}[(2k(r) + rk'(r))\sin^2\theta \\ \qquad + \cos^2\theta]\dot{r}\dot{\varphi} - r\cos\theta\,\dot{r}\dot{\psi} - \frac{r^2}{2}[(k(r)-1)\sin 2\theta]\dot{\theta}\dot{\varphi} + \frac{r^2}{2}\sin\theta\,\dot{\theta}\dot{\psi} = 0, \\ E_5: \ -\frac{r^2}{2}\ddot{\psi} - \frac{r^2}{2}\cos\theta\,\ddot{\varphi} - r\cos\theta\,\dot{r}\dot{\varphi} - r\dot{r}\dot{\psi} + \frac{r^2}{2}\sin\theta\,\dot{\theta}\dot{\varphi} = 0, \end{cases}$$

where $f' = \frac{df}{dr}$, $k' = \frac{dk}{dr}$. Suppose that the vector field $\mathbf{Y} = \overline{\xi}\frac{\partial}{\partial s} + \overline{\eta}^1\frac{\partial}{\partial t} + \overline{\eta}^2\frac{\partial}{\partial r} + \overline{\eta}^3\frac{\partial}{\partial \theta} + \overline{\eta}^4\frac{\partial}{\partial \varphi} + \overline{\eta}^5\frac{\partial}{\partial \psi}$ is the Noether symmetry operator of (4.10) system. Using (4.3) formula, the first prolongation of the vector field $\mathbf{Y}$ will be obtained and then by putting Lagrangian (4.9) and $A = 0$ in equation (4.5) it is found the Noether point symmetries (4.8).

## 5. Conservation laws

When the Noether symmetries of a systems of Euler-Lagrange equations are known, we can obtain the conservation laws for these systems via Noether's theorem [1]. In fact, a conservation law of the system (4.1) is the equation

$$(5.1) \quad D_s T_i = 0$$

on the solutions of (4.1) where $D_s$ is defined in (4.6).

**Theorem 7.** *Suppose that* $\mathbf{Y}$ *is a Noether point symmetry corresponding to a Lagrangian* $\mathcal{L}(s, x^\mu, \dot{x}^\mu)$. *Thus a first integral of (4.7) corresponding to* $\mathbf{Y}$ *is*

$$T = \overline{\xi}\mathcal{L} + (\overline{\eta}^\mu - \dot{x}^\mu \overline{\xi})\frac{\partial \mathcal{L}}{\partial \dot{x}^\mu} - A, \tag{5.2}$$

*where* $A = A(s, x^\mu)$ *is the gauge function.*

Using theorem 7 one will compute all the conserved flows corresponding to the Noether symmetries (4.8). One can get the conserved flow of Noether symmetries (4.8) by following formula

$$T = \overline{\xi}\mathcal{L} + (\overline{\eta}^1 - \dot{t}\overline{\xi})\frac{\partial \mathcal{L}}{\partial \dot{t}} + (\overline{\eta}^2 - \dot{r}\overline{\xi})\frac{\partial \mathcal{L}}{\partial \dot{r}} + (\overline{\eta}^3 - \dot{\theta}\overline{\xi})\frac{\partial \mathcal{L}}{\partial \dot{\theta}} + (\overline{\eta}^4 - \dot{\varphi}\overline{\xi})\frac{\partial \mathcal{L}}{\partial \dot{\varphi}} \tag{5.3}$$
$$+ (\overline{\eta}^5 - \dot{\psi}\overline{\xi})\frac{\partial \mathcal{L}}{\partial \dot{\psi}} - A,$$

where one can put $A = 0$ because the gauge function $A$ is constant. Therefore the conserved flows corresponding to the Noether symmetries (4.8) are

- If $\mathbf{Y}_1 = \dfrac{\partial}{\partial s}$ then

$$T = -f(r)\dot{t}^2 + \frac{k^2(r)}{f(r)}\dot{r}^2 + \frac{r^2}{4}k(r)\dot{\theta}^2 + \frac{r^2}{4}\big[k(r)\sin^2\theta + \cos^2\theta\big]\dot{\varphi}^2 + \frac{r^2}{2}\cos\theta\dot{\varphi}\dot{\psi} + \frac{r^2}{4}\dot{\psi}^2,$$

- If $\mathbf{Y}_2 = \dfrac{\partial}{\partial t}$ then

$$T = -2f(r)\dot{t},$$

- If $\mathbf{Y}_3 = \dfrac{\partial}{\partial \varphi}$ then

$$T = \frac{r^2}{2}\big[k(r)\sin^2\theta + \cos^2\theta\big]\dot{\varphi} + \frac{r^2}{2}\cos\theta\dot{\psi},$$

- If $\mathbf{Y}_4 = \dfrac{\partial}{\partial \psi}$ then

$$T = \frac{r^2}{2}\cos\theta\dot{\varphi} + \frac{r^2}{2}\dot{\psi},$$

- If $\mathbf{Y}_5 = \sin\varphi\dfrac{\partial}{\partial \theta} + \cos\varphi\cot\theta\dfrac{\partial}{\partial \varphi} - \cos\varphi\csc\theta\dfrac{\partial}{\partial \psi}$ then

$$T = \frac{r^2}{2}k(r)\sin\varphi\dot{\theta} + \frac{r^2}{4}\big[k(r) - 1\big]\sin 2\theta\cos\varphi\dot{\varphi} + \frac{r^2}{2}[\cos\theta\cot\theta - \cos\varphi\csc\theta]\dot{\psi},$$

- If $\mathbf{Y}_6 = \cos\varphi\dfrac{\partial}{\partial \theta} - \sin\varphi\cot\theta\dfrac{\partial}{\partial \varphi} + \sin\varphi\csc\theta\dfrac{\partial}{\partial \psi}$ then

$$T = \frac{r^2}{2}k(r)\cos\varphi\dot{\theta} - \frac{r^2}{4}\big[k(r) - 1\big]\sin 2\theta\sin\varphi\dot{\varphi} + \frac{r^2}{2}[\sin\varphi\csc\theta - \cos\theta\sin\varphi\cot\theta]\dot{\psi},$$

where the functions $f(r)$ and $k(r)$ are defined in (1.2).

Department of Mathematics, Faculty of Basic science, Babol University of Technology, Babol, Iran.
*E-mail address*: `r_bakhshandeh@nit.ac.ir`